\documentclass[12pt]{amsart}
\usepackage{amssymb,amsmath,amsthm,latexsym}
\usepackage{mathrsfs}
\usepackage{a4wide}
\usepackage{eucal}
\usepackage[usenames,dvipsnames]{color}
\usepackage[utf8]{inputenc}
\usepackage[T1]{fontenc}
\usepackage{dsfont}

\theoremstyle{remark}
\newtheorem*{maint}{\emph{\textbf{Theorem}}}
\newtheorem*{maintb}{\emph{\textbf{Theorem A}}}

\newtheorem*{remark*}{Remark}
\theoremstyle{definition}

\numberwithin{equation}{section}
\makeatother

\newcommand{\vertiii}[1]{{\left\vert\kern-0.25ex\left\vert\kern-0.25ex\left\vert #1 
    \right\vert\kern-0.25ex\right\vert\kern-0.25ex\right\vert}}

\newcounter{smallromans}

{\end{list}}

\newcounter{smallromansdash}

{\end{list}}

\newcounter{bigromans} 
  {\end{list}}

\title[A letter concerning Leonetti's paper]{A letter concerning Leonetti's paper `Continuous Projections onto Ideal Convergent Sequences'}
\author[T.~Kania]{Tomasz Kania}
\address{Institute of Mathematics, Czech Academy of Sciences, \v{Z}itn\'{a} 25, 115~67 Prague 1, Czech Republic}
\email{tomasz.marcin.kania@gmail.com}

\date{\today}
\thanks{The author acknowledges with thanks funding from funding received from GA\v{C}R project 17-27844S; RVO 67985840 (Czech Republic)}

\begin{document}

\begin{abstract}
    Leonetti proved that whenever $\mathcal I$ is an ideal on $\mathbb N$ such that there exists an~uncountable family of sets that are not in $\mathcal I$ with the property that the intersection of any two distinct members of that family is in $\mathcal I$, then the space $c_{0,\mathcal I}$ of sequences in $\ell_\infty$ that converge to 0 along $\mathcal I$ is not complemented. We provide a shorter proof of a more general fact that the quotient space $\ell_\infty / c_{0,\mathcal I}$ does not even embed into $\ell_\infty$.
\end{abstract}

\maketitle

Very recently, Leonetti (\cite{leonetti}) distilled a property of ideals $\mathcal{I}$ on the set of natural numbers (shared by ideals that are meagre when regarded as subsets of the Cantor set; see \cite[Lemma~2.3]{leonetti}, hence for example by the ideal of sets that have density 0) that gives a fairly satisfactory sufficient condition for non-complementability in $\ell_\infty$, the space of all bounded scalar sequences, of the space $c_{0, \mathcal{I}}$ consisting of sequences that converge to 0 along $\mathcal{I}$. Leonetti's proof is an interesting refinement of Whitley's proof of the Phillips--Sobczyk theorem (\cite{whitley}; see also \cite[Theorem 2.5.4]{ak}), which asserts that the space $c_0$ is not complemented in $\ell_\infty$. (All necessary terminology will be explained in subsequent paragraphs.)\medskip

More specifically, Leonetti's result contributes to the problem of characterisation of those ideals $\mathcal{I}$ of $\mathbb N$ for which the space $c_{0, \mathcal I}$ is complemented in $\ell_\infty$ proposed by P\'erez Hern\'andez during the Winter School in Abstract Analysis 2017 held in Svratka, Czech Republic. Conspicuously, the classical space $c_0$ coincides with $c_{0, \mathcal I}$ where $\mathcal I$ is the ideal of finite sets, so it is not complemented by virtue of the above-mentioned Phillips--Sobczyk theorem. On the other hand, when $\mathcal{I}$ is a maximal ideal (that is when the dual filter is an ultrafilter), the subspace $c_{0, \mathcal{I}}$ has codimension 1 in $\ell_\infty$, so it is complemented. (Similarly, the intersection $\mathcal{I}$ of finitely many maximal ideals will correspond to $c_{0,\mathcal{I}}$ having finite co-dimension in $\ell_\infty$, so in particular $c_{0,\mathcal{I}}$ is complemented in this case.) As every non-principal ideal contains the ideal of finite sets and, at the same time, is contained in a maximal ideal, the problem is quite tantalising indeed. Let us record Leonetti's result formally.\medskip

\begin{maint}Suppose that $\mathcal I$ is an ideal on $\mathbb N$ with the property that there exists an~uncountable family $\mathcal{A}\subset \wp(\mathbb N)\setminus \mathcal{I}$ such that $A \cap B\in \mathcal{I}$ for distinct $A,B\in \mathcal{A}$. Then the space $c_{0, \mathcal I}$ is not complemented in $\ell_\infty$.\end{maint}

We strengthen the above result by noticing that not only is $c_{0, \mathcal I}$ uncomplemented in $\ell_\infty$ but the quotient space $\ell_\infty / c_{0, \mathcal I}$ does not embed into $\ell_\infty$. We have then the following result.

\begin{maintb}Suppose that $\mathcal I$ is an ideal on $\mathbb N$ with the property that there exists an~uncountable family $\mathcal{A}\subset \wp(\mathbb N)\setminus \mathcal{I}$ such that $A \cap B\in \mathcal{I}$ for distinct $A,B\in \mathcal{A}$. Then $\ell_\infty / c_{0, \mathcal I}$ is not isomorphic to a subspace of $\ell_\infty$. In particular, $c_{0, \mathcal I}$ is not complemented in $\ell_\infty$.\end{maintb}

A closed subspace $E$ of a Banach space $X$ is \emph{complemented} whenever there exists a~closed subspace $F$ of $X$ such that $X=E\oplus F$; this, in turn, is equivalent to the existence of a~bounded linear map $T\colon X/E\to X$ such that the composite map $T\pi$ is the identity map when restricted to $E$; here $\pi\colon X\to X/E$ denotes the canonical quotient map.\smallskip

Let $\Gamma$ be a set. A family $\mathcal{I}$ of subsets of $\Gamma$ is an \emph{ideal} on $\Gamma$, when it is closed under finite unions and $B\in \mathcal I$, whenever $B\subseteq A$ and $A\in \mathcal I$. Let $\mathcal I$ be an ideal on the set of natural numbers. A sequence $(x_n)_{n=1}^\infty$ in a metric space $(X,d)$ converges to $x\in X$ along $\mathcal I$, whenever for every $\varepsilon >0$ there is $A\in \mathcal{I}$ such that for every $n\notin A$ we have $d(x_n, x) < \varepsilon$. For every ideal $\mathcal{I}$, the subspace $c_{0, \mathcal I}$ comprising all bounded sequences that converge to 0 along $\mathcal I$ is a closed ideal of $\ell_\infty$, that is a closed subspace which is closed under multiplication by arbitrary elements of $\ell_\infty$. There is a~standard picture of the space $c_{0, \mathcal I}$ as a space of continuous functions on a certain locally compact space that vanish at infinity. \smallskip 

Firstly, let us recall that $\ell_\infty$ is isometrically isomorphic as an algebra to $C(\beta \mathbb N)$, where  $\beta \mathbb N$  is the \v{C}ech--Stone compactification of the integers. Secondly, as $\beta \mathbb N$ consists of all ultrafilters on $\mathbb N$ that is topologised by the base $\{p\in \beta \mathbb N\colon A\in p\}$ ($A\subseteq \mathbb N$), one can consider the open subspace $U_{\mathcal I}$ comprising all ultrafilters that extend the filter dual to $\mathcal{I}$. Set $K_{\mathcal{I}} = \beta \mathbb N \setminus U_{\mathcal {I}}$. By the Tietze--Uryoshn extension theorem, $c_{0, \mathcal {I}}$ is isometric to $C_0(U_\mathcal{I})$, the space of scalar-valued continuous functions on $U_\mathcal{I}$ vanishing at infinity, and $\ell_\infty / c_{0, \mathcal {I}}$ is isometric to $C(K_{\mathcal{I}})$.

\begin{proof}[Proof of Theorem A]Let $\pi\colon \ell_\infty \to \ell_\infty / c_{0,\mathcal{I}}$ be the quotient map. Since $c_{0,\mathcal{I}}$ is an~algebraic ideal of $\ell_\infty$, $\pi$ is, in fact, an algebra homomorphism and $\ell_\infty / c_{0,\mathcal{I}}$ algebraically isomorphic to $C(K_{\mathcal I})$. \smallskip 

For any $A\subseteq \mathbb N$ let us consider the indicator function $\mathds{1}_{A}\in \ell_\infty$. By the hypothesis, $$0=\pi(\mathds{1}_{A\cap B})=\pi(\mathds{1}_{A} \cdot \mathds{1}_{B})=\pi(\mathds 1_A)\cdot \pi (\mathds 1_B)$$ for any distinct $A,B\in \mathcal{A}$, yet $\pi(\mathds{1}_{A}) \neq 0$. Consequently, $\{\pi(\mathds{1}_{A})\colon A\in \mathcal{A}\}$ is a family of pairwise orthogonal non-zero idempotents in $C(K_{\mathcal I})$---that is $\{0,1\}$-valued functions---so it spans an isometric copy of the non-separable space $c_0(\mathcal{A})$ (as $\mathcal{A}$ is uncountable). Consequently, $\ell_\infty / c_{0,\mathcal{I}}$ cannot embed into $\ell_\infty$ (and so $c_{0,\mathcal{I}}$ is not complemented in $\ell_\infty$) as it contains $c_0(\mathcal{A})$. (Indeed, it is a standard fact that $c_0(\Gamma)$ does not embed into $\ell_\infty$ because $\ell_\infty^*\cong \ell_1^{**}$ is weak*-separable by Goldstine's theorem, but $c_0(\Gamma)^*$ is not unless $\Gamma$ is countable and separability is transferred by quotient maps.)\end{proof}

\subsection*{Closing remark} A quick (however unnecessarily high-tech) way of explaining why $c_0$ is not complemented is $\ell_\infty$ is by appealing to the Grothendieck property of the latter space (as proved by Grothendieck himself \cite{Groth}), while observing that $c_0$ clearly lacks it and this property is preserved by surjective linear operators. (A Banach space $X$ is \emph{Grothendieck} whenever every sequence in $X^*$ that converges weak* also converges weakly.) It is then natural to ask the following question.\medskip

Suppose that $\mathcal I$ is an ideal as in the statement of Theorem A. Is it true that $c_{0,\mathcal I}$ is \emph{not} a~Grothendieck space?

\end{document}